\documentclass[12pt]{amsart}
\usepackage{fancyhdr}
\usepackage{amssymb,latexsym,amscd}
\usepackage[cp850]{inputenc}
\newcommand{\cO}{\mathcal{O}}

\newcommand{\lra}{\longrightarrow}

\newcommand{\ra}{\rightarrow}

\newcommand{\CC}{\mathbb C}

\theoremstyle{plain}

\begin{document}
\title[Clifford's Theorem]{On Clifford's Theorem for rank-3 Bundles}
\author{H. Lange}
\author{P.E. Newstead}
\address{Mathematisches Institut\\
              Universit\"at Erlangen-N\"urnberg\\
              Bismarckstra\ss e $1\frac{ 1}{2}$\\
              D-$91054$ Erlangen\\
              Germany}
\email{lange@mi.uni-erlangen.de}
\address{Department of Mathematical Sciences\\University of Liverpool\\
           Peach Street\\Liverpool L69 7ZL\\
           U.K.}
\email{newstead@liv.ac.uk}
\thanks{Supported by DFG Contracts Ba 423/8-1. Both authors are members
of the research group VBAC (Vector Bundles on Algebraic Curves), which is
partially supported by EAGER (EC FP5 Contract no. HPRN-CT-2000-00099) and
by EDGE (EC FP5 Contract no. HPRN-CT-2000-00101).}
\keywords{vector bundle, subbundle}
\subjclass[2000]{Primary:14H60;Secondary:14F05, 32L10}
\begin{abstract}
In this paper we obtain bounds on $h^0(E)$ where $E$ is a semistable bundle of rank 3 over a smooth irreducible 
projective curve $X$ of genus $g \geq 2$ defined over an algebraically closed field of characteristic 0. These
bounds are expressed in terms of the degrees of stability 
$s_1(E), \,\,\, s_2(E)$. We show also that in some cases the bounds
are best possible. These results extend recent work of J. Cilleruelo and I. Sols for bundles of rank 2. 

\end{abstract}
\maketitle


\section{Introduction}

Let $X$ be a smooth irreducible projective curve of genus $g \geq 2$ over an algebraically closed field of 
characteristic 0, and let $E$ be a vector bundle of rank $n$ and degree $d$ over $X$. We recall that $E$ is 
called {\it special} if $h^0(E)$ and $h^1(E)$ are both non-zero.

If $n=1$, the classical Clifford's Theorem provides an upper bound for $h^0(E)$ when $E$ is special, and this 
has been extended to semistable bundles of any rank (see \cite[Theorem 2.1]{2}). Now, for any $E$ of rank $n$
and $1 \leq r \leq n-1$, we define the {\it $r$-th degree of stability} $s_r$ by
$$
s_r(E) = rd - n\max(\deg F),
$$
where the maximum is taken over all subbundles $F$ of rank $r$ of $E$. Note that $E$ is stable (semistable) 
if and only if $s_r(E) > 0$ ($s_r(E) \geq 0$) for all $r$. In case $n=2$, J. Cilleruelo and I. Sols \cite{3} have recently 
obtained a refined version of Clifford's theorem where the bound on $h^0(E)$ depends on $s_1(E)$. Our object 
in this paper is to investigate to what extent this result can be extended to bundles of rank 3.

In section 2, we suppose $E$ is semistable and obtain bounds for $h^0(E)$ 
first in terms of $s_1(F)$,
where $F$ is a rank-2 quotient of minimal degree, and then as a consequence in terms of $s_1(E)$ and
$s_2(E)$ (Theorem 2.3). We then use elementary transformations to show that in some cases our bounds are essentially
best possible. In other cases, the bounds are definitely not best possible, and we give examples also of this situation. 
The necessary properties of elementary transformations are given in section 3, and the 
examples are constructed in section 4. Finally, in section 5, for the 
sake of completeness, we give bounds on 
$h^0(E)$ in the case where $E$ is not semistable.\\

\section{An upper bound for $h^0(E)$ for semistable $E$ }

As in the introduction, let $X$ be a smooth irreducible projective curve of genus $g \geq 2$ defined over an 
algebraically closed field of characteristic 0.
The aim of this section is to give an upper bound on $h^0(E)$ for a semistable vector bundle $E$ of rank 3 
and degree $d$ on $X$ in terms of $d$ and 
the invariants $s_1(E), s_2(E)$ and $s_1(F)$, where $F$ denotes a 
quotient bundle of rank 2 of minimal degree of $E$. For convenience we write
$s_{r} = s_{r}(E)$. \\

{\bf Proposition 2.1}: {\it Let $E$ be semistable of rank $3$ and 
degree $d$ on $X$ and $F$ a 
rank-$2$ quotient of $E$ of minimal degree. 
 If $s_1 \leq 2s_2$ and 
$$
{\max}\left(s_1, \frac{3}{2}s_1(F) - \frac{s_1}{2}\right) \leq d \leq 6g - 6 - \frac{3}{2}s_1 (F) - \frac{s_1}{2},  \eqno(2.1)
$$
 then}
$$
h^0(E) \leq \left[\frac{d}{2} - \frac{s_1(F)}{2}\right] + 3.
$$

{\it Proof}: Write $L = Ker(E \lra F)$ and let $M$ be a line 
subbundle of $F$ of maximal degree. Then $\deg L = \frac{d-s_1}{3}$ and 
$\deg M = \frac{d}{3} + \frac{s_1}{6} - \frac{s_1(F)}{2}$. This gives a commutative diagram 
$$
\begin{array}{ccccccccc}
0 & \lra & L & \lra & N & \lra & M & \lra & 0\\
&&\| && \downarrow && \downarrow &&\\
0 & \lra & L & \lra & E & \lra & F & \lra & 0
\end{array}
$$
where the upper row is the pullback of the lower row via the 
inclusion $M \hookrightarrow F$.
We claim that $F$ is semistable.
 
To see this, note that $\deg N = \frac{1}{6} (4d - s_1 - 3s_1(F))$. 
But $N$ is a rank-2 subbundle of $E$. Hence $\deg N \leq \frac{2d - s_2}{3}$
by definition of $s_2$. This implies
$$
s_1(F) \geq \frac{2s_2 - s_1}{3} \geq 0,       \eqno(2.2)
$$ 
that is $F$ is semistable. The assumption (2.1) implies  
$$s_1(F) \leq \deg(F) \leq 4g - 4 - s_1(F).$$ 
Hence \cite[Theorem 0.2]{3} applies to give
$$
h^0(F) \leq \frac{\deg F - s_1(F)}{2} + 2.        \eqno(2.3)
$$ 
On the other hand 
$$
0 \leq \frac{d-s_1}{3} = \deg(L) \leq 2g - 2 - \frac{1}{2}(s_1(F) + s_1).
$$ 
Since $E$ and $F$ are semistable, we have $s_1(F) + s_1 \geq 0$.
Thus $L$ is in the range of Clifford's Theorem for line bundles 
and we obtain

\begin{eqnarray*}
h^0(E) \leq h^0(L) + h^0(F) &\leq &\frac{d - s_1}{6} + 1 + \frac{d}{3} + \frac{s_1}{6} - \frac{s_1(F)}{2} + 2\\  
& \leq & \frac{d}{2} - \frac{s_1(F)}{2} + 3.
\end{eqnarray*}   

\begin{flushright}
$\square$
\end{flushright}

{\bf Remark 2.2}: The estimate in Proposition 2.1 can be slightly improved by using instead of (2.3) 
the full version of the theorem of Cilleruelo-Sols (see \cite[Theorem 0.2]{3})
as follows.

The {\it Krawtchouk polynomials} $K_r(n,N)$ are defined by the identity
$$
\sum_{r=0}^{\infty} K_r(n,N)z^r = (1-z)^n(1+z)^{N-n}.
$$  
Then we have with the assumptions of Proposition 2.1
$$
h^0(E) \leq \frac{d}{2} - \frac{s_1(F)}{2} + 2 + \delta     \eqno(2.4)$$
where
$$ \delta = \left\{ \begin{array}{lll}
0 & if & K_{\frac{1}{6}(2d+s_1-3s_1(F))+1}(g,2g-s_1(F)) \not= 0\\
1 & if & K_{\frac{1}{6}(2d+s_1-3s_1(F))+1}(g,2g-s_1(F)) = 0.
\end{array} \right.
$$
For non-hyperelliptic curves one would get an even better inequality by using the full version of Clifford's 
Theorem for line bundles or even the Clifford index.
On the other hand we will see in section 4 that the inequality (2.4) is sharp for $X$ hyperelliptic.\\

\vspace{0.3cm}
If $d< s_1$, then $h^0(E) = 0$ by definition of $s_1$. 
If $d > 6g-6 - s_2$, then $h^1(E) = h^0(E^* \otimes \omega_X) = 0$ since $\deg(E^* \otimes \omega_X) =
-d + 6g - 6 < s_2 = s_1(E^* \otimes \omega_X)$, 
and thus $h^0(E)$ may be computed by Riemann-Roch. Hence we may assume 
$$
s_1 \leq d \leq 6g - 6 - s_2.              \eqno(2.5)
$$ 
The following theorem gives an upper bound for $h^0(E)$ for every 
semistable $E$ with $\deg E = d$ in this range.\\   

{\bf Theorem 2.3}: {\it Let $E$ be semistable of rank $3$ and degree $d$.}
\renewcommand{\theenumi}{\roman{enumi}}
\renewcommand{\labelenumi}{(\theenumi)}
\begin{enumerate}
\item {\it If
\begin{enumerate}
\item $\frac{s_1}{2} \leq s_2 \leq 2s_1$ and $s_1 \leq d \leq 6g-6-s_2$ or 
\item $s_2 > 2s_1$ and $ s_2 - s_1 \leq d \leq 6g-6-s_2$ or
\item $s_2<\frac{s_1}{2}$ and $s_1\leq d\leq 6g-6-(s_1-s_2),$
\end{enumerate} 
then} 
$$
 \qquad \qquad h^0(E) \leq \left[\frac{d}{2} - \frac{1}{6} \max(2s_2 - s_1,2s_1 - s_2)\right] + 3
$$
\item {\it If $s_2 > 2s_1$ and $s_1 \leq d < s_2 - s_1$, then }
$$
h^0(E) \leq \left[\frac{d}{2} - \frac{s_1}{2}\right] + 1,
$$
\item {\it If $ s_2 < \frac{s_1}{2}$ and $6g-6-(s_1 - s_2) < d \leq 6g-6-s_2$, then} 
$$
h^0(E) \leq \left[\frac{d}{2} - \frac{s_2}{2}\right] + 1.
$$
\end{enumerate}

{\it Proof}: As in Proposition 2.1 let $L$ be a line subbundle of $E$ of maximal degree with quotient $F$. 
The proof of (i) proceeds in several steps.\\

{\bf Step 1}: {\it If $s_1 \leq 2s_2$ and $\max(s_1, s_2-s_1) \leq d \leq 6g-6-s_2$, then}
$$
h^0(E) \leq \frac{d}{2} - \frac{2s_2 - s_1}{6} + 3.      \eqno(2.6)
$$
For the proof suppose first that in addition to $s_1 \leq 2s_2$ we have also (2.1). 
Then according to Proposition 2.1 and (2.2) 
$$
h^0(E) \leq \frac{d}{2} -\frac{s_1(F)}{2} + 3 \leq \frac{d}{2} - \frac{2s_2 - s_1}{6} +3.
$$
If $s_2 -s_1 \leq d < \frac{3}{2}s_1(F) - \frac{s_1}{2}$ then $\deg F < s_1(F)$ 
and thus $h^0(F) = 0$. Hence
$$
h^0(E) = h^0(L) \leq \frac{d-s_1}{6}+1 = \frac{d}{2} - \frac{2d+s_1}{6} + 1 \leq \frac{d}{2} - \frac{2s_2 - s_1}{6} + 1.
$$
If $6g-6 -\frac{3}{2}s_1(F) - \frac{s_1}{2} < d \leq 6g - 6 - s_2$ 
then $h^1(F) = 0$ and so by Riemann-Roch $h^0(F) = \frac{2d+s_1}{3} +2(1-g)$ implying
\begin{eqnarray*}
h^0(E) \leq h^0(L) +h^0(F) & \leq & \frac{d-s_1}{6} + 1 + \frac{2d+s_1}{3} + 2(1-g)\\
& = & \frac{d}{2} + \frac{s_1}{6} + \frac{d-6g+6}{3} + 1\\ 
& \leq & \frac{d}{2} - \frac{2s_2 - s_1}{6} + 1.
\end{eqnarray*}
This completes the proof of Step 1.\\

{\bf Step 2}: {\it If $s_2 \leq 2s_1$ and $s_1 \leq d \leq 6g - 6 - \max(s_2, s_1-s_2)$, then}
$$
h^0(E) \leq \frac{d}{2} - \frac{2s_1 - s_2}{6} + 3.      \eqno(2.7)
$$
For the proof note that passing from $E$ to $E^* \otimes \omega_X$ 
interchanges $s_1$ and $s_2$. So Step 1 gives 
$h^0(E^* \otimes \omega_X) \leq - \frac{d}{2} + 3g - 3 - 
\frac{2s_1 - s_2}{6} + 3$ and thus

\begin{eqnarray*}
h^0(E) & = & h^0(E^* \otimes \omega_X) + d + 3 - 3g\\
& \leq & \frac{d}{2} - \frac{2s_1 - s_2}{6} + 3.
\end{eqnarray*}

{\bf Step 3}: 
\renewcommand{\theenumi}{\alph{enumi}}
\begin{enumerate}
\item If $ \frac{s_1}{2} \leq s_2 \leq 2 s_1$, both formulas (2.6) and (2.7) apply in the full range $s_1 \leq d \leq 6g-6 - s_2$.\\ 
\item If $s_2 > 2s_1$, formula (2.6) applies in the range $s_2 - s_1 \leq d \leq 6g - 6 - s_2$, since then $s_1 < s_2 - s_1$. 
But in this case $\max(2s_2 - s_1, 2s_1 - s_2) = 2s_2 - s_1$.\\
\item Finally, if $s_2 < \frac{s_1}{2}$, formula (2.7) applies in the range $s_1 \leq d \leq 6g-6-(s_1 - s_2)$, 
since then $s_1 - s_2 > s_2$. But in
this case $\max(2s_2 - s_1, 2s_1 - s_2) = 2s_1 - s_2$. 
\end{enumerate}

This completes the proof of (i).

{\it Proof of} (ii): Suppose $s_2 > 2s_1$ and  $d < s_2 - s_1$, then by (2.2)
$s_1(F) \geq \frac{2s_2-s_1}{3} > \frac{2d+s_1}{3} = \deg F$
implying $h^0(F) = 0$. Since $d \geq s_1$ we get

\begin{eqnarray*}
h^0(E) = h^0(L) & \leq & \frac{d-s_1}{6} + 1\\
& \leq &\frac{d}{2} - \frac{s_1}{2} + 1
\end{eqnarray*}
for all $d$ in the range $s_1 \leq d < s_2 - s_1$.

{\it Proof of} (iii): This is exactly similar to the proof of (ii).  
\hfill  $\square$\\

{\bf Remark 2.4}: If $E$ admits a rank-2 quotient $F$ of minimal degree such that the Krawtchouk polynomial satisfies
$K_{\frac{1}{6}(2d+s_1-3s_1(F))+1}(g,2g-s_1(F)) \not= 0$, then Remark 2.2 can be applied to give a slightly better bound in 
case (i) of Theorem 2.3:
$$
h^0(E) \leq \left[\frac{d}{2} - \frac{1}{6} \max(2s_2 - s_1,2s_1 - s_2)\right] + 2.\eqno(2.8)
$$\\

\section{Elementary transformations}

In order to construct some vector bundles with a large space of global sections we need some 
properties of elementary transformations which we collect in this section. We state them for bundles of arbitrary rank, although we 
need them here only for bundles of rank 3.

Let $E$ be a vector bundle of rank $n$ and degree $d$ over the curve $X$. 
For any point $x \in X$ we denote by $E(x)$ the fibre of $E$ and
by $\CC_x$ the skyscraper sheaf on $X$ with fibre $\CC$ at $x$ and 0 elsewhere. By an {\it elementary transformation} 
of $E$ we shall mean  
a vector bundle $E'$ fitting into an exact sequence
$$
0 \lra E \lra E' \lra \CC_x \lra 0.      \eqno(3.1)
$$
The elementary transformation $E'$ of $E$ determines a 1-dimensional subspace $l_x$ of $E(x)$, namely the kernel of the induced map
$E(x) \lra E'(x)$. Conversely, any 1-dimensional subspace $l_x \subset E(x)$ determines an elementary transformation $E'$ of $E$
as follows: Let $H_x$ denote the hyperplane of the dual vector space $E^*(x)$ defined by $l_x$ and $\CC_x$ the skyscraper sheaf 
with fibre 
$E^*(x)/H_x$ at $x$. Let $E'{}^*$ denote the kernel of the canonical map 
$E^* \lra \CC_x$. Its dual $E' = E'{}^{**}$ fits 
into an exact sequence (3.1). We call the vector bundle $E'$ the elementary transformation of $E$ 
{\it associated to} $l_x \subset E(x)$.

For any $r$, $1 \leq r \leq n-1$, and any vector bundle $E$ of rank $n$ and degree $d$ we denote by $d^r_{max}(E)$ the maximal degree 
of a subbundle of rank $r$ of $E$, that is 
$$
d^r_{max}(E) = \frac{1}{n}(rd - s_r(E)).
$$
Moreover we define for any integer $i \geq 0$
$$
SB^{\,i}_r(E) = \{ \mbox{subbundles} \,\, F \,\, \mbox{of} \,\,  E \,\, \mbox{of rank}\,\, r \,\, 
\mbox{and degree} \,\, d^r_{max}(E) - i \}.
$$
If $E$ is of degree $d$, this can be given a scheme structure by identifying it with 
an open set of Grothendieck's scheme $Quot_{n-r,d - d^r_{max}(E) + i}(E)$ 
of quotients of $E$ of rank $n-r$ and degree $d - d^r_{max}(E) + i$. 
Note that $SB^{\,0}_r(E)$ is the set of rank-$r$ subbundles of maximal degree of $E$
which is a projective scheme.\\

The proof of the 
following lemma is straightforward (see e.g. \cite[Proposition 1.6]{1}).\\

{\bf Lemma 3.1}: {\it Let $E'$ denote the elementary transformation of $E$ 
associated to $l_x \subset E(x)$. Then for} $1 \leq r \leq n-1$
$$
s_r(E') = \left\{ \begin{array}{l} s_r(E) - (n-r)\\
                                       s_r(E) + r \end{array} \right.
                                       if \quad
 \begin{array}{l} \exists F \subset SB^{\,0}_r(E) \,\, with \,\, l_x \subset F(x),\\
                         l_x \not\subset F(x) \,\, for \,\, all \,\, F \subset SB^{\,0}_r(E). \end{array}
$$   \\

Let $E'$ be any elementary transformation of $E$ with exact sequence (3.1). It is easy to see that 
the set of subbundles $F'$ of $E'$ is in canonical bijection to 
the set of subbundles $F$ of $E$ via the map $F' \lra  F = F' \cap E$. 
If $F$ is a subbundle of $E$, we always denote by $F'$ the corresponding
subbundle of $E'$. With this notation the sets $SB^0_r(E)$ and $SB^0_r(E')$ are related 
as in the following lemma the proof of which is straightforward.\\

{\bf Lemma 3.2} (i): {\it If $s_r(E') = s_r(E) - (n-r)$, then}
$$
SB^{\,0}_r(E') = \{F' \subset E' | F \in SB^{\,0}_r(E) \,\, and \,\, l_x \subset F(x)\}
$$
(ii): {\it if $s_r(E') = s_r(E) + r$, then}
$$
SB^{\,0}_r(E') = \{F' \subset E' | F \in SB^{\,0}_r(E) \} \, \cup \, \{F' \subset E' |  F \in SB^{\,1}_r(E) 
\, and \, l_x \subset F(x) \}.
$$\\

An elementary transformation of $E$ is defined by a pair $(x,l_x)$ where $x$ 
is a point of $X$ and $l_x$ a line (= 1-dimensional subspace) of
the vector space $E(x)$. Hence the set of elementary transformations $elm(E)$ 
forms a projective bundle of fibre dimension $r-1$ over the curve $X$. In particular 
$elm(E)$ is a projective variety of dimension $r$ and it makes sense to speak of a {\it general} elementary transformation.\\

{\bf Proposition 3.3}: {\it If $\dim SB^0_r(E) < n-r$, then for a general elementary transformation $E'$ of $E$}
$$
s_r(E') = s_r(E) + r.
$$

{\it Proof}: Using the identification of $SB^0_r(E)$ with a Quot-scheme, we see that there exists a bundle $U$ over
$X\times SB^0_r(E)$ which is universal as a family of subbundles of $E$ of rank $r$ and degree $d^r_{max}(E)$.
Hence we get a canonical morphism
$$
\varphi: U \lra E
$$
where $U$ and $E$ are considered as varieties. Now
$$
r + 1 + \dim SB^0_r(E) = \dim U < r + 1 + (n - r) = \dim E. 
$$
Hence a general point $p \in E$ is not contained in the image of $\varphi$. 
If $p \in E(x)$ and $l_x$ denotes the line in $E(x)$ spanned by $p$, the
elementary transformation $E'$ of $E$ associated to $l_x$ satisfies 
$s_r(E') = s_r(E) + r$ according to Lemma 3.1.\hfill $\square$\\     

Now write $E = E_0$. We will inductively construct sequences of elementary transformations 
$$
0 \lra E_{k} \lra E_{k+1} \lra \CC_{x_k} \lra 0    \eqno(3.2)
$$
associated to $l_{x_k} \subset E_k(x_k)$ for all $k \geq 0$. As above  
the set of subbundles of rank $r$ of $E_0$ is in canonical bijection to 
the set of subbundles of rank $r$ of $E_k$ via these exact sequences.
If $F_0$ is a subbundle of $E_0$ of rank $r$, we denote by $F_k$ the corresponding subbundle of $E_k$ for all $k$. \\

{\bf Proposition 3.4}: {\it Suppose that for some positive integer $m$  the vector bundle $E_0$ satisfies
$$
\dim SB^{\,i}_r(E_0) < (i+1)(n-r)
$$
for $i=0,\ldots,m-1$. Then for a general sequence of 
elementary transformations (3.2) we have}
$$
s_r(E_m) = s_r(E_0) + mr.
$$\\

The proof follows immediately by induction from the following lemma.\\

{\bf Lemma 3.5}: {\it Suppose that for $0 \leq k \,\,(\leq m-1)$ the sequence of
elementary transformations (3.2) is  constructed up to $E_k$ in 
such a way that 
$$
s_r(E_k) = s_r(E_0) + kr \,\, \mbox{and} \,\, \dim SB^{\,i}_r(E_k) < (i+1)(n-r)
$$
for $i = 0, \ldots , m - 1 -k$. Then we have for the elementary transformation $E_{k+1}$ of $E_{k}$ 
associated to a general $x_k \in X$ and a general line
$l_{x_{k}} \subset E_{k}(x_{k})$
$$
s_r(E_{k+1}) = s_r(E_0) +(k+1) r \,\, \mbox{and} \,\, \dim SB^i_r(E_{k+1}) < (i+1)(n-r)
$$ 
for} $i=0, \ldots, m - 1 - (k+1)$.\\

{\it Proof}: Since $\dim SB^0_r(E_k) < n-r$, we have by Proposition 3.3 for a general $x_k \in X$ 
and a general line $l_{x_k} \subset E(x_k)$
$$
s_r(E_{k+1}) = s_r(E_k) + r = s_r(E_0) + (k+1)r.
$$
Moreover a slight generalization of Lemma 3.2 implies
$$
\begin{array}{lll}
SB^{\,i}_r(E_{k+1}) &= &\{F_{k+1} \subset E_{k+1} | F_k \in SB^{\,i}_r(E_k) \quad and \quad l_{x_k} \not\subset F_k(x) \}\\ 
&& \, \cup \, \{F_{k+1} \subset E_{k+1} | F_k \in SB^{\,i+1}_r(E_k) \quad and \quad l_{x_k} \subset F_k(x) \}.
\end{array}
$$
Now 
$$
\dim\{F_{k+1} \subset E_{k+1} | F_k \in SB^{\,i}_r(E_k) \, and \, l_{x_k} \not\subset F_k(x) \} \le \dim SB^i_r(E_k) < (i+1)(n-r)
$$
for $i = 0, \ldots , m - 1 -k$ and 
$$
\dim \{F_{k+1} \subset E_{k+1} | F_k \in SB^{\,i+1}_r(E_k)\} = \dim SB^{\,i+1}_r(E_k) < (i+2)(n-r)
$$ 
for $i = 0, \ldots, m-1 -(k+1)$. But for a subbundle $F_k \subset E_k$ of rank $r$ to contain the line $l_{x_k}$ 
imposes $n-r$ conditions. Hence
$$
\dim \{F_{k+1} \subset E_{k+1} | F_k \in SB^{\,i+1}_r(E_k) \quad and \quad l_{x_k} \subset F_k(x) \} < (i+2)(n-r) -(n-r)
$$
for $i = 0, \ldots, m-1 -(k+1)$.\hfill $\square$\\

Now let $n$ be equal to 2 or 3. We want to apply Proposition 3.4 in order 
to construct some vector bundles of rank $n$. Let $p_1, \ldots, p_n$ be
general points of the curve $X$ and consider the vector bundle
$$
E_0 = \cO_X(p_1) \oplus \cdots \oplus \cO_X(p_n).
$$
$E_0$ is of rank and degree $n$ with $d^1_{max}(E_0) = 1$.\\

{\bf Lemma 3.6}: $\dim SB_1^i(E_0) < (i+1)(n-1)$ {\it for} $n=2, 3$ {\it and }
$i = 0, \ldots, g-1$.\\  

{\it Proof}: We may assume $i > 0$, the assertion for $i=0$ being obvious. 
Let 
$$Quot_{n+i-1}(E_0)$$ denote the Quot-scheme 
parametrizing quotients of rank $n-1$ and degree $n+i-1$ of $E_0$. 
Since $SB_1^i(E_0)
\subset Quot_{n+i-1}(E_0)$, it suffices to show that 
$$\dim Quot_{n+i-1}(E_0)<(i+1)(n-1)\eqno(3.3)$$ locally at all points of 
$SB^i_1(E_0)$.

Suppose $L \in SB_1^i(E_0)$ with associated exact sequence
$$
0 \lra L \lra E_0 \lra F \lra 0.
$$
We have to estimate the dimension of $H^0(F \otimes L^{-1})$, 
since this is the tangent space of $Quot_{n+i-1}(E_0)$ at the point $L$.

Suppose first $n=2$. Then $F \otimes L^{-1} \simeq L^{-2}(p_1 + p_2)$ 
is a line bundle of degree $2i \leq 2g - 2$. Hence by Clifford's Theorem for 
line bundles
$$h^0(L^{-2}(p_1 + p_2)) \leq \frac{deg(L^{-2}(p_1 + p_2))}{2} = i,$$ 
unless either $i=g-1$ and $L^{-2}(p_1+p_2)\simeq \omega_X$ or 
$X$ is hyperelliptic and $L^{-2}(p_1 +p_2) \simeq h^i$ where $h$ denotes 
the hyperelliptic line bundle. For all other $L$
the local dimension of $SB_1^i(E_0)$ at $L$ is at most  $i$.
On the other hand, there are only finitely many exceptional cases, and in all
these cases
$$h^0(E_0 \otimes L^{-1}) = h^0(L^{-1}(p_1)) + h^0(L^{-1}(p_2)) 
\leq \frac{i+1}{2} + \frac{i+1}{2} = i+1.$$
When $X$ is not hyperelliptic and $L^{-2}(p_1+p_2)\simeq \omega_X$, this
is clear. When $X$ is hyperelliptic, note that $L^{-1}(p_1)$ and
$L^{-1}(p_2)$ cannot both be powers of $h$ since $p_1$, $p_2$ are general 
points of $X$; hence neither can be a power of $h$. So the subset of
$SB^i_1(E_0)$ given by inclusions of $L$ in $E_0$ has dimension $\leq i$. 
This proves (3.3) for $n=2$.

Finally suppose $n=3$. 
Then $F$ is a rank-2 vector bundle of degree $i+2$. It is easy to see that $s_1(F) \geq -i$. 
On the other hand $s_1(F) \leq i$, since 
at least one of the line subbundles of degree 1 of $E_0$ maps to a nonzero 
subsheaf of $F$.
Let $M$ be a maximal subbundle of $F \otimes L^{-1}$. Since $\deg (F \otimes L^{-1}) = 3i$, we have 
$$
i \leq \deg M \leq 2i.
$$
Since also $0 < 2i \leq 2g-2$, both $M$ and $N = (F \otimes L^{-1})/M$ lie in the range of Clifford's Theorem for 
line bundles. So 
$$
\begin{array}{ll}
h^0(F \otimes L^{-1}) & \leq h^0(M) + h^0(N)\\
& \leq [\frac{3i}{2}] + 2 \leq 2i + 1.
\end{array} 
$$
This proves (3.3) for $n=3$ and completes the proof of the lemma.
\hfill$\square$

\medskip
{\bf Corollary 3.7}: {\it Suppose $n = 2$ or $n=3$ and $E_m$ is obtained from
$E_0$ by a general sequence of elementary transformations  (3.2) 
with $1\le m\le g$. Then} $s_1(E_m) = m$.\\

{\it Proof}: Since $s_1(E_0) = 0$, this follows at once from Lemma 3.6 
and Proposition 3.4.

\hfill $\square$\\

\section{Examples}

In this section we give some examples showing that some of the estimates of Section 2 are  sharp. For this 
assume that $X$ is a hyperelliptic 
curve of genus $g$ and denote by $h$ the unique line bundle on $X$ of 
degree 2 with $h^0(h) = 2$. (If $X$ is not hyperelliptic, then the bounds 
of Section 2 cannot be attained for $h^0(E)$ as noted at the end of 
Remark 2.2.) Note that for a hyperelliptic curve it follows from 
\cite[Proposition 3]{4} (see also the two paragraphs following the
statement of Theorem 0.2 in \cite{3})
that, if $F$ is a bundle of rank $2$ with
$0<s_1(F)\le\deg F\le 4g-4-s_1(F)$, then

$$h^0(F)\le\frac{\deg F-s_1(F)}{2}+1.\eqno(4.1)$$ 
Hence, by the proof of Proposition 2.1, if $s_1(E) \leq 2s_2(E)$, (2.1) holds, and $s_1(F)>0$, then
$$h^0(E)\le\left[\frac{d}2-\frac{s_1(F)}2\right]+2.\eqno(4.2)$$
(This is just (2.4) with $\delta=0$.) Moreover, a careful analysis of the
proof shows that, under the hypotheses of Theorem 2.3(i), and provided
$s_1$ and $s_2$ are not both zero, then (2.8) holds. Note that these
improvements are independent of the values of the Krawtchouk polynomials.

\medskip  
\begin{center}
{\bf (a) Examples with $s_1 = s_2 = 0$ }
\end{center}  

Start with $F_0 = \cO(p_1) \oplus \cO(p_2)$ for $p_1, p_2 \in X$ and let $F_m$ be the bundle obtained from $F_0$ by a general sequence 
of elementary transformations (3.2) for some $m \leq g$. Then we have an exact sequence $0 \lra F_0 \lra F_m \lra T_m \lra 0$ 
with a torsion sheaf $T_m$ of length $m$.

Suppose now $m = 4n + 2$ and consider for $k = 0, \dots, g - 2 - \frac{m}{2}$ the vector bundle
$$
F_{m,k} = F_m \otimes h^k.
$$
Then $s_1(F_{m,k}) = m$ by Corollary 3.7, and $\deg F_{m,k} = m+ 4k + 2 = 4(n+k+1) \leq 4g - 4 - m$. Hence, by (4.1),
$$h^0(F_{m,k}) \leq 2k+2.$$
On the other hand,
tensoring the above exact sequence by $h^k$, we get
$$
0 \lra h^k(p_1) \oplus h^k(p_2) \lra F_{m,k} \lra T_m \lra 0.
$$
This implies $h^0(F_{m,k}) \geq h^0(h^k(p_1) \oplus  h^k(p_2)) = 2k + 2$ and thus
$$
h^0(F_{m,k}) = 2k+2. 
$$

Now let $L = h^{n+k+1}$ and consider 
$$
E = L \oplus F_{m,k}.
$$
Then $\deg E = 6(n+k+1)$ and $s_1(E) = s_2(E) =0$. Moreover 
$$
h^0(E) = n+k+2 + 2k+2 = n + 3k + 4 = \frac{\deg E}{2} -
\frac{s_1(F_{m,k})}{2} + 2,
$$
which is just the bound (4.2). Hence the estimate (4.2) is best possible in this case.\\

{\bf Remark 4.1}: If $g=2$, there are no allowable values of $k$, so this method gives examples only for $g\ge3$.\\

\begin{center}
{\bf (b) Stable examples}
\end{center}  

Start with $E_0 = \cO(p_1) \oplus \cO(p_2) \oplus \cO(p_3)$ for $p_1, p_2, p_3 \in X$ and let $E_m$ be the bundle obtained from $E_0$ 
by a general sequence of elementary transformations (3.2) for some $m \leq g$. Then we have the following diagram
$$
\begin{array}{ccccccccc}
&&0&&0&&&&\\
&&\downarrow&&\downarrow&&&&\\
0& \lra & \cO(p_1)& \lra &\cO(p_1)& \lra&0&&\\
&&\downarrow&&\downarrow&&\downarrow&&\\ 
0&\lra&E_0&\lra&E_m&\lra&T_m&\ra& 0\\
&&\downarrow&&\downarrow&&\|&&\\ 
0&\lra&\cO(p_2) \oplus \cO(p_3)&\lra&F&\lra&T_m&\lra&0.\\
&&\downarrow&&\downarrow&&\downarrow&&\\ 
&&0&&0&&0&&
\end{array} 
$$
Here $T_m$ is a torsion sheaf of length $m$.
By Corollary 3.7, $s_1(E_m) = m$ and  $\cO(p_1)$ is a subbundle of maximal degree of $E_m$.
Moreover, in the lower exact sequence the vector bundle $F$ is obtained by a general sequence of 
elementary transformations starting from $\cO(p_2) \oplus \cO(p_3)$.
In particular $s_1(F) = m$ according to Corollary 3.7. 
Hence the middle horizontal exact sequence of the diagram and inequality 
(4.2) imply
$$
3 = h^0(E_0) \leq  h^0(E_m)
 \leq  \left[\frac{\deg E_m}{2} - \frac{s_1(F)}{2}\right] + 2=3.
$$
Hence if the hypotheses of (4.2) are satisfied, we conclude that $h^0(E_m)=3$ and
the bound (4.2) is sharp in this case.\\ 

Thus in order to get many stable examples, it remains to show the following lemma,
since inequality (2.1)  holds provided $g\ge3$.\\

{\bf Lemma 4.2}: {\it Suppose $m$ is even and $2 \leq m \leq g$. Then $s_1(E_m) \leq 2s_2(E_m)$.}\\

{\it Proof}: According to Corollary 3.7 $s_1(E_m) = m$. So it suffices to show that $s_2(E_m) \geq \frac{m}{2}$.
For this it is enough to show that $s_2(E_m) \geq \frac{m-3}{2}$ for all $m$, since then $s_2(E_m) \equiv 2m \bmod 3$ 
implies that $s_2(E_m) \geq \frac{m}{2}$ for $m$ even.

First we claim that, if $s_2(E_m) < \frac{m}{2}$, then $E_m$ admits only finitely many maximal rank-2 subbundles. 
For the proof let $G$ denote a maximal rank-2 subbundle of $E_m$ and $L = E_m/G$. It suffices to show that $H^0(G^* \otimes L) = 0$, since 
this is the tangent space of the corresponding Quot-scheme. The assumption implies that $\deg(G^* \otimes L) < \frac{m}{2}$. So if 
$H^0(G^* \otimes L) \not= 0$,
we would get 
$$s_1(G) = s_1(G^* \otimes L) \leq \deg(G^* \otimes L) < \frac{m}{2}.$$
But a maximal line subbundle $M$ of $G$ is also a subbundle of $E_m$. So $s_1(E_m) = m$ gives $\deg M \leq 1$ and thus
$s_1(G) > \frac{m}{2}$, a contradiction. 

Returning now to the proof of the assertion $s_2(E_m) \geq \frac{m-3}{2}$, note first that it
is certainly valid for $E_0$. Suppose it holds for $E_m$. If $s_2(E_m) < \frac{m}{2}$, then $E_m$ admits only finitely many 
maximal rank-2 subbundles as we have seen above. Hence, by Proposition 3.3, 
$$
s_2(E_{m+1}) = s_2(E_m) + 2 \geq \frac{(m+1) - 3}{2}.
$$ 
If $s_2(E_m) \geq \frac{m}{2}$, then 
$$s_2(E_{m+1}) \geq s_2(E_m) - 1 \geq \frac{(m+1) -3}{2}.
$$ This completes the proof by induction on $m$. \hfill $\square$\\

{\bf Remark 4.3}: 
One can show that Lemma 4.2 is true also for $m= 1,3,5$ and it is possible that it is true for all $m$.
Note also that, if $m=g=2$, (2.1) fails. On the other hand, $g=2$, $m=1$ is allowed, so we do get some examples even for $g=2$.\\

\begin{center}
{\bf (c) Examples for Theorem 2.3}
\end{center}

Constructing examples to illustrate Theorem 2.3 (or rather (2.8)) seems to be harder. However, we can at least construct
a few examples. As above, let $X$ be a hyperelliptic curve and consider the 
vector bundle 
$$
E_0 = \cO(p_1) \oplus \cO(p_2) \oplus \cO(p_3)
$$
with general points $p_1, p_2, p_3 \in X$. Let $E_1$ 
be the general elementary transformation of $E_0$. Since $E_0$ has just 3 maximal subbundles of rank $2$, it is an immediate
consequence of Proposition 3.3 that $s_2(E_1)=2$. Of course, we already know that $s_1(E_1)=1$ from Corollary 3.7.

Now define for
$k=0,\ldots,g-2$
$$E_{1,k}=E_1\otimes h^k.$$
We have $\deg E_{1,k}=6k+4$. Moreover
$$h^0(E_{1,k})\ge h^0(h^k(p_1))+h^0(h^k(p_2))+h^0(h^k(p_3))=3k+3.$$
On the other hand $s_1\le\deg E_{1,k}\le 6g-8=6g-6-s_2$. 
Hence (2.8) yields 
$$
h^0(E_{1,k}) \leq 3k + 3,
$$
and the bound of Remark 2.4 is sharp in this case.\\

Now let $E_2$ be obtained from a general sequence of elementary transformations (3.2) with $m=2$.
For
$k=0,\ldots,g-2$, define
$$E_{2,k}=E_2\otimes h^k.$$
Then $\deg E_{2,k}=6k+5$ and a calculation similar to the above gives
$$h^0(E_{2,k})\ge3k+3.$$
In fact, since $h^1(E_{1,k})>0$ and $E_{2,k}$ is a general elementary transformation of $E_{1,k}$, $h^0(E_{2,k})=3k+3$. 
In this case we know that $s_1(E_{2,k})=2$, but\\

{\bf Lemma 4.4}: $s_2(E_2) = 1$.\\

{\it Proof}: We prove first that 
$$\dim SB_2^1(E_0)=2.\eqno(4.3)$$
For this, we must show that the family of rank-2 subbundles of $E_0$ of degree $1$ has dimension $2$. It is simpler to
work with quotient line bundles of degree 2. For the line bundle $h$, note first that $h$ can arise as a quotient of $E_0$. 
Then, by a simple calculation,
$$h^0(E_0^*\otimes h)=3,$$
from which it follows that $h$ gives rise to a 2-dimensional family as required. On the other hand, if $L\not\simeq h$ is
a quotient of $E_0$ of degree $2$, 
it is easy to show that  $L\simeq{\cO}(p_i+p_j)$ for some $i\ne j$ and then $h^0(E_0^*\otimes L)=2$. Since
there are finitely many such $L$, this completes the proof of (4.3). 

For a general choice of $l_x$, the condition that $F\in SB_2^1(E_0)$ contains $l_x$ imposes just one condition on $F$.
It follows from (4.3) and Lemma 3.2(ii) that $E_1$ has a $1$-dimensional family of maximal subbundles of rank $2$. But then Lemma 3.1
implies that
$$s_2(E_2)=s_2(E_1)-1=1.$$
\hfill $\square$\\

We still have $s_1\le\deg E_{2,k}\le 6g-6-s_2$, so (2.8) gives $h^0(E_{2,k})\le3k+4$. Thus $E_{2,k}$ does not give the exact bound.

In fact, when $k=0$, we have $d=5$, $s_1=2$ and $s_2=1$, and in this case the bound of (2.8) cannot be attained for $g\ge3$ (even
if $E$ is not obtained in the above manner). Indeed $E$ is necessarily stable and has slope less than $2$. So by \cite{5} it follows
that
$$h^0(E)\le 3+\frac1g(d-3)=3+\frac2g.$$
For $g=2$, on the other hand, this bound can be attained; again by \cite{5} there exists a stable $E$ of degree $5$ with 
$h^0(E)=4$. Since $s_1\equiv5\bmod3$ and $s_2\equiv10\bmod3$, it is easy to see that the only possible values of $s_1$, $s_2$ are 
$s_1=2$, $s_2=1$. So $E$ attains the bound (2.8).\\

\section{Upper bound for $h^0(E)$ for unstable $E$ }

For the sake of completeness we give in this section also an upper bound 
for $h^0(E)$ in the case of an 
unstable (i.e. not semistable) vector bundle $E$ of rank $3$ and 
degree $d$ on the curve $X$, now no longer assumed to be hyperelliptic.   
Recall that $E$ is unstable if $s_1(E) < 0$ or $s_2(E) < 0$. We may assume 
$$
s_1 = s_1(E) < 0.
$$
Indeed, if we have a bound in this case, we also get 
a bound in case $s_1(E) \geq 0$, since then $s_1(E^* \otimes \omega_X) = s_2(E) < 0$ and
$$
h^0(E) = h^0(E^* \otimes \omega_X) + d + 3 -3g.  
$$
Hence a bound for $h^0(E^* \otimes \omega_X)$ gives a bound for $h^0(E)$.\\

According to the assumptions there is an exact sequence
$$
0 \lra L \lra E \lra F \lra 0
$$
with a line bundle $L$ of degree $\frac{d-s_1}{3} > \frac{d}{3}$ and 
a rank-2 vector bundle $F$ of degree $\frac{2d+s_1}{3}$. 
As in the proof of Proposition 2.1,
we see that
$$
s_1(F) \geq \frac{2s_2 - s_1}{3}.      \eqno(5.1)
$$\\

Suppose first that $F$ is semistable. As in (2.5) we may assume $s_1 \leq d \leq 6g - 6 - s_2$. 
With these assumptions we have the following result:\\

{\bf Proposition 5.1}: {\it If $s_1 < 0$ and $F$ is semistable, then} 
$$
h^0(E) \leq h^0(L) + h^0(F)
$$
{\it with}
$$
h^0(L)  \left\{ \begin{array}{l}
\leq \vspace{0.1cm} \frac{1}{6}(d-s_1) + 1 \\ = \frac{1}{3}(d-s_1) +1-g \end{array} \,\, {\it if} \,\,
\begin{array}{l} \vspace{0.1cm} s_1 \leq d \leq 6g-6+s_1\\
 6g-6+s_1 < d \leq 6g-6-s_2
\end{array} \right.
$${\it and}
$$
h^0(F)
\left\{ \begin{array}{l}
= \vspace{0.1cm} 0\\
\leq \vspace{0.1cm} \frac{1}{3}(d+s_1-s_2) + 2\\
= \vspace{0.1cm} \frac{1}{3}(2d+s_1) + 2 - 2g
\end{array} \right. 
\,\, if \,\,
\begin{array}{l}
\vspace{0.1cm} s_1 \leq d < \frac{3}{2}s_1(F) -\frac{s_1}{2}\\
\vspace{0.1cm} \frac{3}{2}s_1(F) -\frac{s_1}{2} \leq d \leq 6g-6-\frac{3}{2}s_1(F) -\frac{s_1}{2}\\
6g-6-\frac{3}{2}s_1(F) -\frac{s_1}{2} < d \leq 6g-6-s_2.
\end{array}
$$\\

{\it Proof}: Note first that $\deg L = \frac{d-s_1}{3} \ge 0$, so either 
Clifford's Theorem or Riemann- Roch gives
the estimate for $h^0(L)$.
On the other hand $\frac{3}{2}s_1(F) - \frac{s_1}{2} \leq d \leq 6g-6- \frac{3}{2}s_1(F) - \frac{s_1}{2}$ is equivalent to 
$s_1(F) \leq deg(F) \leq 4g-4-s_1(F)$. Hence \cite[Theorem 0.2]{3} 
implies in this case
$$
h^0(F) \leq \frac{1}{3}(d+s_1-s_2) + 2.
$$
If $d < \frac{3}{2}s_1(F) - \frac{s_1}{2}$, then $\deg F<s_1(F)$ and so 
$h^0(F) = 0$. If $6g-6-\frac{3}{2}s_1(F) -\frac{s_1}{2} < d \leq 6g-6-s_2$, 
then we may apply Riemann-Roch to give $h^0(F) = \frac{1}{3}(2d +s_1) +2 -2g$. 
Now $h^0(E) \leq h^0(L) + h^0(F)$ gives the assertion.\hfill $\square$\\

Again one obtains a slightly better estimate for $h^0(E)$ by applying the full version of \cite[Theorem 0.2]{3} 
using the Krawtchouk polynomials.
It is easy to see that these bounds are best possible by considering suitable 
direct sums $E = L \oplus F$ with $h^0(L)$ and $h^0(F)$ maximal.\\

Finally let us assume that $F$ is unstable. Using Clifford's Theorem and Riemann-Roch we obtain the following result:\\

{\bf Proposition 5.2}: If $s_1 < 0$ and $F$ is unstable, then
$$
 h^0(E) \leq h^0(L) + h^0(F)
$$   {\it with}\\
$$
h^0(L)  \left\{ \begin{array}{l}
\leq \vspace{0.1cm} \frac{1}{6}(d-s_1) + 1 \\ = \frac{1}{3}(d-s_1) +1-g \end{array} \,\, {\it if} \,\,
\begin{array}{l} \vspace{0.1cm} s_1 \leq d \leq 6g-6+s_1\\
 6g-6+s_1 < d \leq 6g-6-s_2
\end{array} \right.
$$  \,\, {\it and}
$$
h^0(F)  \left\{ \begin{array}{l}
 = \vspace{0.1cm} 0\\
 \leq \vspace{0.1cm} \frac{1}{6}(d+s_1-s_2) + 1\\
 \leq \vspace{0.1cm} \frac{1}{6}(2d+s_1) + 2 \\
\leq \vspace{0.1cm} \frac{1}{12}(6d+4s_1-2s_2)-g+2\\
= \frac{1}{3}(2d +s_1) +2 - 2g
\end{array}
\,\,if\,\,
\begin{array}{l}
\vspace{0.2cm} s_1 \leq d < \frac{3s_1(F) -s_1}{2}\\
\vspace{0.1cm}\frac{3s_1(F)-s_1}{2} \leq d < -\frac{3s_1(F) + s_1}{2}\\
\vspace{0.1cm} -\frac{3s_1(F) + s_1}{2} \leq d \leq 6g-6 +\frac{3s_1(F) - s_1}{2}\\
\vspace{0.1cm} 6g-6 +\frac{3s_1(F) - s_1}{2} < d \leq 6g-6 -\frac{3s_1(F) + s_1}{2}\\
6g-6 -\frac{3s_1(F) + s_1}{2} < d \leq 6g - 6 - s_2.
\end{array} \right.
$$\\

Again it is easy to see using direct sums of suitable line bundles on hyperelliptic curves that the bounds in Proposition 5.2 are sharp.

\end{document}